\documentclass[12pt]{article}
\usepackage{latexsym}

\usepackage[tbtags]{amsmath}
\usepackage{epsfig,amstext,amssymb,amsthm,latexsym}
\usepackage{amssymb,color}

\definecolor{c20}{rgb}{0.,0.7,0.}
\definecolor{c30}{rgb}{0.,0.,1.}
\definecolor{c40}{rgb}{1,0.1,0.7}
\definecolor{c50}{rgb}{1,0,0}

\date{}
\newtheorem{theorem}{Theorem}[section]
\newtheorem{lemma}{Lemma}[section]
\newtheorem{proposition}{Proposition}[section]
\newtheorem{corollary}{Corollary}[section]
\newtheorem{remark}{Remark}[section]
\newtheorem{example}{Example}[section]

\makeatletter 
\@addtoreset{equation}{section}
\makeatother 
 
\begin{document}

\title{On the Borel-Cantelli lemma}

\author{  Alexei Stepanov
\thanks{Department of Mathematics, Izmir University of Economics,
35330, Balcova, Izmir, Turkey; {\it  alexeistep45@mail.ru;
alexei.stepanov@ieu.edu.tr} }, \small{\it Izmir University of
Economics, Turkey}}

\def\abstractname{}

\date{\begin{abstract} In the present note,
we generalize  the first part of  the Borel-Cantelli lemma. By this generalization, we obtain some strong limit results.
\end{abstract}}

\maketitle  \vspace{3mm}\noindent


\section{Introduction}
Suppose $A_1,A_2,\cdots$ is a sequence of events on a common
probability space and that $A^c_i$ denotes the complement of event
$A_i$. The Borel-Cantelli lemma, presented below as Lemma~\ref{lemma1.1},
is used extensively for producing strong limit theorems.
\begin{lemma}\label{lemma1.1}
\begin{enumerate}
\item If, for any sequence  $A_1,A_2,\cdots$ of events,
\begin{equation}\label{1.1}
\sum_{n=1}^\infty P(A_n)<\infty,
\end{equation}
then $P(A_n\ i.o.)=0$, where i.o. is an abbreviation for
"infinitively often``.
\item If $A_1,A_2,\cdots$ is a sequence of
independent events and if
\begin{equation}\label{1.2}
\sum_{n=1}^\infty P(A_n)=\infty,
\end{equation}
then $P(A_n\ i.o.)=1$.
\end{enumerate}
\end{lemma}
The independence condition in the second part of the
Borel-Cantelli lemma is weakened by a number of authors, including
Chung and Erdos (1952), Erdos and Renyi (1959), Lamperti (1963),
Kochen and Stone (1964), Spitzer (1964), Ortega and Wschebor
(1983), and Petrov (2002), (2004). One can also refer to
Martikainen and Petrov (1990), and Petrov (1995) for related
topics.

The first part of the Borel-Cantelli lemma is generalized in
Barndorff-Nielsen (1961) (Lemma~\ref{lemma1.2}) and Balakrishnan and Stepanov (2010) (Lemma~\ref{lemma1.3}).
\begin{lemma}\label{lemma1.2}
Let $A_1,A_2,\ldots$ be a sequence of events such that
$P(A_n)\rightarrow 0$. If
\begin{equation}\label{1.3}
\sum_{n=1}^\infty P(A_n A^c_{n+1})<\infty,
\end{equation}
then $P(A_n\ i.o.)=0$.
\end{lemma}
\begin{lemma}\label{lemma1.3}
Let  $A_1,A_2,\ldots$ be a sequence of events such that
$P(A_n)\rightarrow 0$.  If
\begin{equation}\label{1.4}
\sum_{n=1}^\infty P(A^c_n A_{n+1})<\infty,
\end{equation}
then $P\{A_n\ i.o.\}=0$.
\end{lemma}

In the present work, by simple rewriting conditions in (\ref{1.3}) and (\ref{1.4}) we obtain a new generalization of the first part of the Borel-Cantelli lemma.
We  will  show that $P(A_n\ i.o.)=0$ might be  also true when $P(A_n)\rightarrow 0$ and inequality (\ref{1.2}) holds. This generalization is illustrated by some limit results.

The rest of this paper is organized as follows. In Section~2, we present our  results.
In Section~3, we use the results of Section~2 for deriving  strong limit theorems for dependent random variables.

\section{Results}

\begin{lemma}\label{lemma2.1}
Let $A_1, A_2, \dots$ be a sequence of events such that
$P(A_n)\rightarrow 0$.  Let (\ref{1.2}) hold true,
\begin{equation}\label{2.1}
\sum_{n=1}^\infty P(A_n A_{n+1})=\infty
\end{equation}
and
\begin{equation}\label{2.2}
\sum_{n=1}^\infty\left[P(A_n)-P(A_n A_{n+1})\right]<\infty.
\end{equation}
Then $P(A_n\ i.o.)=0$.
\end{lemma}
\begin{gproof}{of Lemma~\ref{lemma2.1}} Observe that $P(A_nA^c_{n+1})=P(A_n)-P(A_nA_{n+1})$. The truth of Lemma~\ref{lemma2.1}  follows from Lemma~\ref{lemma1.2}.
\end{gproof}
\begin{remark}\label{remark2.1} Condition (\ref{2.2}) in Lemma~\ref{lemma2.1} might be replaced by the condition $$
\sum_{n=1}^\infty\left[P(A_{n+1})-P(A_n A_{n+1})\right]<\infty.
$$
In that case we could prove Lemma~\ref{lemma2.1} by using Lemma~\ref{lemma1.3} instead of Lemma~\ref{lemma1.2}.
\end{remark}
\section{Applications: limit results}
The importance of  the theoretical results of Section~2 is shown in this section.
The following limit theorem can be easily derived  from Lemma~\ref{lemma2.1}.
\begin{theorem}\label{theorem3.1} Let $X_1,X_2,\ldots$ be a sequence of dependent random variables such that  $X_n\stackrel{p}{\rightarrow}\mu $, where $\mu$ is a constant. Let  for all small $\varepsilon >0$,
$$
\sum_{n=1}^\infty P(X_n\not\in[\mu -\varepsilon,\mu +\varepsilon ])=\infty,
$$
$$
\sum_{n=1}^\infty P(X_n\not\in[\mu -\varepsilon,\mu +\varepsilon ], X_{n+1}\not\in[\mu -\varepsilon,\mu +\varepsilon ])=\infty,
$$
and
$$
\sum_{n=1}^\infty \left[P(X_n\not\in[\mu -\varepsilon,\mu +\varepsilon ])-P(X_n\not\in[\mu -\varepsilon,\mu +\varepsilon ], X_{n+1}\not\in[\mu -\varepsilon,\mu +\varepsilon ])\right]<\infty.
$$
Then
$$
X_n\stackrel{a.s.}{\rightarrow}\mu.
$$
\end{theorem}
Theorem~\ref{theorem3.1} is illustrated by  Example~\ref{example3.1} below. Here, we also formulate two trivial statements: Proposition~\ref{proposition3.1} and Corollary~\ref{corollary3.1}. These statements are given to simplify the presentation of Example~\ref{example3.1}.
\begin{proposition}\label{proposition3.1}
Let $A_1,A_2,\cdots$ be a decreasing sequence of  events such that  $P(A_n)\rightarrow0$. Then $P(A_n\ i.o.)=0$.
\end{proposition}
Proposition~\ref{proposition3.1} can  be easily proved directly. Corollary~\ref{corollary3.1} follows from Proposition~\ref{proposition3.1}.
\begin{corollary}\label{corollary3.1}
Let $X_1,X_2,\cdots$ be an ordered sequence of  random variables such that  $X_n\stackrel{p}{\rightarrow}\mu $, where $\mu $ is a constant. Then $X_n\stackrel{a.s.}{\rightarrow}\mu$.
\end{corollary}
\begin{example}\label{example3.1} Let $X_1,\dots,X_n,\ldots$ be a sequence of dependent random variables defined for any $n\geq 1$ by the Clayton copula
$$
F(x_1,x_2,\ldots,x_n)=
\left[\frac{1}{x_1}+\frac{1}{x_2}+\dots+\frac{1}{x_n}-(n-1)\right]^{-1},
$$
where  $0<x_i<1\ (1\leq i\leq n)$.  Let $M_n=\max\{X_1,\ldots,X_n\}$. Then, for any $x\in(0,1)$ and $n\rightarrow \infty $
$$
P(M_n\leq x)=\left[n\left(\frac{1}{x}-1\right)+1\right]^{-1}\rightarrow 0.
$$
We see that $M_n$ converges in probability to 1. Since $\sum_{n=1}^\infty P(M_n\leq x)=\infty $, we can not apply the first part of the Borel-Cantelli lemma for deriving the corresponding strong limit result. However, $M_n$ is an ordered sequence of  random variables. By Corollary~\ref{corollary3.1}, we have
$$
M_n\stackrel{a.s.}{\rightarrow}1.
$$

Theorem~\ref{theorem3.1} allows us to obtain a more elaborate strong limit result. We will show that
\begin{equation}\label{3.2}
M_n^{n^\alpha}\stackrel{a.s.}{\rightarrow}1\quad (0<\alpha <1).
\end{equation}
Observe first that
$$
P(M_n^{n^\alpha }\leq x)=\left[n\left(x^{-n^{-\alpha }}-1\right)+1\right]^{-1}\sim
$$$$
(-\log x)n^{\alpha -1}\rightarrow 0.
$$
It follows that $M_n^{n^\alpha}\stackrel{p}{\rightarrow}1$. However, we can not
apply here Proposition~\ref{proposition3.1}, since $M_n^{n^\alpha}$ is not an ordered sequence. Instead, we utilize  Theorem~\ref{theorem3.1}. Observe that
$$
\sum_{n=1}^\infty P(M_n^{n^\alpha }\leq x)=\infty,
$$
$$
\sum_{n=1}^\infty P(M_n^{n^\alpha }\leq x,M_{n+1}^{(n+1)^\alpha }\leq x)=\sum_{n=1}^\infty\frac{1}{n(x^{-n^{-\alpha}}-1)+x^{-(n+1)^{-\alpha}}}=\infty
$$
and
$$
\sum_{n=1}^\infty [P(M_n^{n^\alpha }\leq x)-P(M_n^{n^\alpha }\leq x, M_{n+1}^{(n+1)^\alpha }\leq x)]=
$$
$$
\sum_{n=1}^\infty \frac{x^{-(n+1)^{-\alpha}}-1}{ (n(x^{-n^{-\alpha}}-1)+1) (n(x^{-n^{-\alpha}}-1)+x^{-(n+1)^{-\alpha}})}\sim
$$
$$
\frac{1}{(-\log x)}\sum_{n=1}^\infty \frac{1}{n^{2-\alpha}}<\infty,
$$
The convergence in (\ref{3.2}) readily follows.
\end{example}

\section*{References}
\begin{description} {\small
\item Balakrishnan, N., Stepanov, A. (2010).\ Generalization of
the Borel-Cantelli lemma. {\it The Mathematical Scientist}, {\bf
35} (1), 61--62.

\item Barndorff-Nielsen, O. (1961).\ On the rate of growth of the
partial maxima of a sequence of independent identically
distributed random variables. {\it Math. Scand.},  {\bf 9},
383--394.

\item Chung, K.L. and Erdos, P. (1952).\ On the application of the
Borel-Cantelli lemma. {\it Trans.  Amer. Math. Soc.},  {\bf 72},
179--186.

\item Erdos, P. and Renyi, A. (1959).\ On Cantor's series with
convergent $\sum 1/q_n$. {\it Ann. Univ. Sci. Budapest. Sect.
Math.}, {\bf 2}, 93--109.

\item Kochen, S.B. and Stone, C.J. (1964).\ A note on the
Borel-Cantelli lemma. {\it Illinois J. Math.},  {\bf 8}, 248--251.

\item Lamperti, J. (1963).\ Wiener's test and Markov chains. {\it
J. Math. Anal. Appl.},  {\bf 6}, 58--66.

\item Martikainen, A.I., Petrov, V.V., (1990).\ On the
Borel–Cantelli lemma. {\it Zapiski Nauch. Semin. Leningrad. Otd.
Steklov Mat. Inst.}, {\bf 184}, 200–-207 (in Russian). English
translation in:  (1994). {\it J. Math. Sci.}, {\bf 63}, 540–-544.

\item Petrov, V.V. (1995).\ {\it Limit Theorems of Probability
Theory}. Oxford University Press, Oxford.

\item Petrov, V.V. (2002).\ A note on the Borel-Cantelli lemma.
{\it Statist. Probab. Lett.}, {\bf 58}, 283--286.

\item Petrov, V.V. (2004).\ A generalization of the Borel-Cantelli
Lemma. {\it Statist. Probab. Lett.}, {\bf 67}, 233--239.

\item Ortega, J., Wschebor, M., (1983).\ On the sequence of
partial maxima of some random sequences. {\it Stochastic Process.
Appl.}, {\bf 16}, 85–98.

\item Spitzer, F. (1964).\ {\it Principles of Random Walk}.  Van
Nostrand, Princeton, New Jersey.

}
\end{description}
\end{document}